\documentclass[letterpaper, 10 pt, conference]{ieeeconf}  

\IEEEoverridecommandlockouts                              

\usepackage{amsmath} 
\usepackage{amssymb}  
\usepackage{mathtools}
\let\labelindent\relax
\usepackage{multirow}
\usepackage{xcolor}
\usepackage{subcaption}
\usepackage{float}
\usepackage{balance}
\usepackage[shortlabels]{enumitem}
\newcommand*{\QEDB}{\hfill\ensuremath{\square}}%

\newtheorem{asmp}{Assumption}
\newtheorem{rem}{Remark}
\newtheorem{definition}{Definition}

\newtheorem{theorem}{Theorem}

\title{\LARGE \bf
On Event-Triggered Extremum Seeking via Standard and Lie-Bracket Averaging: A Hybrid Dynamical Systems Approach
}

\author{Mahmoud Abdelgalil, Jorge I. Poveda
\thanks{*This work was not presented at any conference. This work was supported in part by the grants NSF ECCS CAREER 2305756 and AFOSR YIP: FA9550-22-1-0211.}
\thanks{The authors are with the Department of Electrical and Computer Engineering,
        University of California, San Diego, 92093, USA
        {\tt\small }}%
}

\begin{document}

\maketitle
\thispagestyle{empty}
\pagestyle{empty}


\begin{abstract}
  We introduce and analyze the stability of a class of event-triggered extremum-seeking algorithms designed to solve resource-aware, model-free, optimization problems. Leveraging recent advances in Lie-Bracket Averaging for hybrid systems, we demonstrate that the proposed controllers can be formulated as well-posed multi-time-scale hybrid systems that satisfy key regularity, stability, and robustness properties. In extremum-seeking systems, exploration and exploitation are inherently coupled. This coupling necessitates careful consideration in the design of the event-triggered controller. To address this challenge, we incorporate a low-pass filter into the algorithm and carefully design the flow and jump sets of the resulting hybrid system. The resulting controller renders the optimal point semi-globally practically asymptotically stable with solutions exhibiting a uniform semi-global dwell time. We also demonstrate how the proposed event-triggered scheme can be modified to allow analysis using traditional averaging tools for hybrid systems by introducing two independent tunable parameters in the controller. Numerical simulations are presented to validate and illustrate the theoretical results.

\end{abstract}

\section{INTRODUCTION}
Extremum seeking (ES) control ranks among the most prominent model-free online optimization and stabilization approaches introduced in the last century \cite{Leblanc1922}. Its ease of implementation, theoretical guarantees of stability, and robustness have made ES systems a compelling choice for a broad array of practical control problems \cite{ariyur2003real,scheinker2024100,yilmaz2024prescribed} in different areas such as source seeking for mobile robots \cite{poveda2021robust,suttner2023nonlocal,abdelgalil2023singularly,bajpai2024model}, control of energy systems \cite{ghaffari2014power,ebegbulem2017distributed}, real-time optimization of engines \cite{popovic2006extremum}, and decision-making in games \cite{frihauf2011nash,krilavsevic2021learning,poveda2022fixed}. Moreover, recent works have illustrated potential connections between seeking dynamics and certain biological systems \cite{abdelgalil2022sea,eisa2023analyzing}. For a recent survey, we refer the reader to \cite{scheinker2024100}.

In an independent development, event-triggered control \cite{heemels2012introduction,postoyan2014framework} has been established as a principled framework to optimize the use of computational and communication resources in feedback control systems. Event-triggered controllers are inherently hybrid, as they combine continuous-time and discrete-time dynamics. Over the past few decades, event-triggered control has been extensively studied, with successful applications in regulation, tracking, optimization, and estimation.

Nevertheless, event-triggered extremum-seeking systems have received comparatively less attention. ES controllers form a class of adaptive algorithms that aim to locate, in real time, the optimal operating point of a system using only output measurements and a high-frequency periodic dithering signal applied to the input. Event-triggered ES algorithms have only recently been explored in \cite{rodrigues2025event} using sampled data structures. However, while general frameworks for the design of hybrid ES based on classical averaging and Lie-Bracket Averaging have been presented in \cite{poveda2017framework} and \cite{abdelgalil2025lie}, respectively, the full potential of event-triggered ES designs remains largely untapped.

Motivated by this gap, in this paper, we introduce a new class of event-triggered ES algorithms formulated within the hybrid dynamical systems framework \cite{goebel2012hybrid}, and grounded in Lie-Bracket Averaging theory \cite{durr2013examples,abdelgalil2025lie}, as well as in traditional averaging theory for hybrid extremum seeking \cite{poveda2017framework}. Specifically, building on a recent result for Lie-Bracket Averaging in hybrid inclusions \cite{abdelgalil2025lie}, we show that the behavior of the proposed controller can be rigorously analyzed via the associated hybrid Lie-bracket averaged dynamics. The properties of closeness of solutions and semi-global practical stability are studied using tools from hybrid inclusions and are shown to extend naturally to the proposed controllers. These properties rely on working with solutions defined on hybrid time domains, for which suitable (graphical) distance metrics can be leveraged for the purpose of analysis. Our main result establishes conditions under which the event-triggered ES controller accurately emulates its hybrid averaged dynamics and achieves semi-global practical asymptotic stability of the desired equilibrium point or optimal set. Zeno behavior is also ruled out by design of the target averaged system. Finally, we show that the proposed scheme can be modified to obtain event-triggered ES controllers suitable for analysis via traditional averaging methods for hybrid systems and hybrid extremum seeking control \cite{poveda2017framework}, thus providing a unifying framework for the study of event-triggered ES systems, illustrated via two particular designs.

The rest of the paper is organized as follows. In Section II we present preliminaries on notation and hybrid systems. In Section III, we study the target hybrid system that we seek to approximate with the Lie-bracket event-triggered controller. Our main control design and results are presented in Section IV. Section V presents numerical simulations, and finally, Section VI concludes the paper.

\section{Preliminaries}
\subsection{Notation}
The bilinear form $\langle \cdot,\cdot\rangle:\mathbb{R}^n\times\mathbb{R}^n\rightarrow\mathbb{R}$ is the canonical inner product on $\mathbb{R}^n$, i.e., $\langle x, y\rangle = x^\top y$, for any two vectors $x,y\in\mathbb{R}^n$. We use $\|x\|=\sqrt{\langle x,x\rangle}$ to denote the $2$-norm of a vector $x\in\mathbb{R}^n$. We denote by $\mathbb{R}_{+}$ the set of positive real numbers, and by $\mathbb{R}_{\geq 0}$ the set of non-negative real numbers. For a compact set $\mathcal{A}\subset\mathbb{R}^n$ and a vector $x\in\mathbb{R}^n$, we use $|x|_{\mathcal{A}}:=\min_{\tilde{x}\in\mathcal{A}}\|x-\tilde{x}\|$. For a continuous function $V:\mathbb{R}^n\rightarrow\mathbb{R}$, we use $\mathcal{L}_c(V)$ to denote the level set $\{x\in\mathbb{R}^n~|~V(x)=c\}$ for any $c\in\mathbb{R}$. A function $\beta:\mathbb{R}_{\geq0}\times\mathbb{R}_{\geq0}\to\mathbb{R}_{\geq0}$ is of class $\mathcal{KL}$ if it is strictly increasing in its first argument, strictly decreasing in its second argument, $\lim_{r\to0^+}\beta(r,s)=0$ for each $s\in\mathbb{R}_{\geq0}$, and  $\lim_{s\to\infty}\beta(r,s)=0$ for each $r\in\mathbb{R}_{\geq0}$. Throughout the paper, for two (or more) vectors $u,v \in \mathbb{R}^{n}$, we write $(u,v)=[u^{\top},v^{\top}]^{\top}$. A function $f$ is said to be $\mathcal{C}^k$ if all of its derivatives up to the $k$th-derivative exist and are continuous.
\subsection{Hybrid Dynamical Systems}
The algorithms studied in this paper integrate both continuous-time and discrete-time dynamics. As such, they will be modeled as hybrid dynamical systems (HDS) aligned with the framework of \cite{goebel2012hybrid}, and characterized by the following equations:
\begin{align}\label{eq:hds_notation}
\mathcal{H}:~~\begin{cases}
    ~~\xi\in C, & \dot{\xi}\hphantom{^+}=F(\xi),\\
    ~~\xi\in D, & \xi^+=G(\xi),
\end{cases}
\end{align}
where $F:\mathbb{R}^n\rightarrow\mathbb{R}^n$ is called the flow map, $G:\mathbb{R}^n\rightarrow\mathbb{R}^n$ is called the jump map, $C\subset\mathbb{R}^n$ is called the flow set, and $D\subset\mathbb{R}^n$ is called the jump set. We use $\mathcal{H}=(C,F,D,G)$ to denote the \emph{data} of the HDS and impose the following standard assumption on $\mathcal{H}$:
\begin{asmp}\label{asmp:regularity}
The sets $C$ and $D$ are closed, and the functions $F$ and $G$ are continuous. \hfill $\QEDB$
\end{asmp}

Any HDS that satisfies Assumption \ref{asmp:regularity} is a \emph{well-posed} HDS in the sense of \cite[Definition 6.29]{goebel2012hybrid}, see \cite[Theorem 6.30]{goebel2012hybrid}. 

\subsection{Solutions to HDS and Zeno behavior}
Solutions to  \eqref{eq:hds_notation} are parameterized by a continuous-time index $t\in\mathbb{R}_{\geq0}$, which increases continuously during flows, and a discrete-time index $j\in\mathbb{Z}_{\geq0}$, which increases by one during jumps. Therefore, solutions to \eqref{eq:hds_notation} are defined on \emph{hybrid time domains} (HTDs). A set $E\subset\mathbb{R}_{\geq0}\times\mathbb{Z}_{\geq0}$ is called a \textsl{compact} HTD if $E=\cup_{j=0}^{J-1}([t_j,t_{j+1}],j)$ for some finite sequence of times $0=t_0\leq t_1\leq \ldots\leq t_{J}$. The set $E$ is an HTD if for all $(T,J)\in E$, $E\cap([0,T]\times\{0,\ldots,J\})$ is a compact HTD. A hybrid arc $\xi$ is a function defined on an HTD. In particular, $\xi:\mathrm{dom}(\xi)\to \mathbb{R}^n$ is such that $\xi(\cdot, j)$ is locally absolutely continuous for each $j$ such that the interval $I_j:=\{t:(t,j)\in \mathrm{dom}(\xi)\}$ has a nonempty interior. A hybrid arc $\xi:\mathrm{dom}(\xi)\to \mathbb{R}^n$ is a solution $\xi$ to the HDS \eqref{eq:hds_notation} if $\xi(0, 0)\in C\cup D$, and:
\begin{enumerate}
\item For all $j\in\mathbb{Z}_{\geq0}$ such that $I_j$ has nonempty interior: $\xi(t,j)\in C$ for all $t\in I_j$, and $\dot{\xi}(t,j)= F(\xi(t,j))$ for all $t\in I_j$.
\item For all $(t,j)\in\mathrm{dom}(\xi)$ such that $(t,j+1)\in \mathrm{dom}(\xi)$: $\xi(t,j)\in D$ and $\xi(t,j+1)= G(\xi(t,j))$.
\end{enumerate}
A solution $\xi$ is said to be \emph{maximal} if it cannot be further extended. A solution $\xi$ is said to be \emph{complete} if the length of its domain is infinite, i.e., for every $T>0$ there exists $(t,j)\in\mathrm{dom}(\xi)$ such that $t+j>T$.

In the analysis of event-triggered systems, it is important to characterize the presence or absence of Zeno solutions---those for which the time domain is bounded in the \( t \)-direction but unbounded in the \( j \)-direction, with increasingly shorter flow intervals between jumps \cite[Ch.2]{goebel2012hybrid}. Indeed, as shown in \cite{postoyan2014framework}, for several common event-triggered schemes it is difficult to rule out Zeno behavior in the entire space, or to guarantee a constant dwell time between jumps in the solutions. Therefore, the following definition, borrowed from \cite[Def. 3]{postoyan2014framework}, will be relevant in the study of Zeno behavior in the event-triggered ES systems considered in this paper.

\vspace{0.1cm}
\begin{definition}
The solutions to \eqref{eq:hds_notation} are said to have a \emph{uniform semi-global dwell time} if for any $\Delta\geq0$, there exists $d(\Delta)>0$ such that for any solution $\xi$ with $|\xi(0,0)|\leq\Delta$ and for any $(s,i),(t,j)\in\mathrm{dom}(\xi)$ with $s+i\leq t+j$, the following holds:
\begin{equation}\label{dwelltime}
j-i\leq\frac{t-s}{d(\Delta)}+1.
\end{equation}
\hfill $\QEDB$
\end{definition}

\subsection{Closeness of Solutions in HDS}
In the analysis and design of extremum-seeking systems using averaging tools, a key mathematical property that enables the transfer of stability results from the averaged system to the actual controller is the so-called ``closeness'' of solutions between the two systems. This property is typically established in the limit as certain tunable parameters are decreased. However, since solutions to hybrid dynamical systems (HDS) are generally discontinuous, standard notions of distance are not well-suited for establishing this closeness; see \cite[Example 3]{goebel2009hybrid} for a counterexample. Instead, to determine whether two solutions to \eqref{eq:hds_notation} are ``close,'' we analyze the distance between their trajectories using appropriate metrics. This is formalized in the following definition:

\vspace{0.1cm}
\begin{definition}
given $T\geq 0$ and $\epsilon>0$, two hybrid arcs $x$ and $y$ are said to be $(T,\epsilon)$-close if the following conditions are satisfied:
\begin{enumerate}
\item For each $(t,j)\in\mathrm{dom}(x)$ with $t+j\leq T$, there exists $s\in\mathbb{R}_{\geq0}$ such that $(s,j)\in\mathrm{dom}(y)$, $|t-s|\leq \epsilon$, and
\begin{equation}
|x(t,j)-y(s,j)|<\epsilon.
\end{equation}
\item For each $(t,j)\in\mathrm{dom}(y)$ with $t+j\leq T$, there exists $s\in\mathbb{R}_{\geq0}$ such that $(s,j)\in\mathrm{dom}(x)$, $|t-s|\leq \epsilon$, and 
\begin{equation}
|y(t,j)-x(s,j)|<\epsilon.
\end{equation}
\end{enumerate}
\hfill $\QEDB$
\end{definition}
One of the main results of this paper is to show that each solution of the proposed event-triggered ES controller is $(T,\epsilon)$-close to some solution of its average target hybrid dynamics (c.f., Theorems 2 and 3) by leveraging Assumption 1 and the analytical results of \cite{abdelgalil2025lie} and \cite{poveda2017framework}.

\subsection{Boundedness and Stability Notions for HDS}
The following stability notions for $\mathcal{H}$, stated with respect to sets, will be relevant in the study of the event-triggered algorithms considered herein.

\vspace{0.1cm}
\begin{definition}
    A compact set $\mathcal{A}\subset C\cup D$ is said to be Uniformly Globally Asymptotically Stable (UGAS) for the HDS $\mathcal{H}$ if there exists a class-$\mathcal{KL}$ function $\beta$ such that every solution $\xi$ of the HDS $\mathcal{H}$ satisfies
    \begin{align*}
        |\xi(t,j)|_{\mathcal{A}}\leq \beta(|\xi(0,0)|_{\mathcal{A}},t+j),
    \end{align*}
    for all $(t,j)\in\mathrm{dom}(\xi)$.  \hfill $\QEDB$
\end{definition}

\vspace{0.1cm}
Finally, when the data of the HDS $\mathcal{H}$ depend on a small parameter $\varepsilon$, we will need the following stability notion.

\vspace{0.1cm}
\begin{definition}
    A compact set $\mathcal{A}\subset C\cup D$ is said to be Semi-Globally Practically Asymptotically Stable (SGPAS) as $\varepsilon\rightarrow 0^+$ for the HDS $\mathcal{H}$ if there exists a class-$\mathcal{KL}$ function $\beta$ such that for every compact set $K$ and every $\kappa > 0$, there exists $\varepsilon^*>0$ such that, for all $\varepsilon\in(0,\varepsilon^*)$, every solution $\xi$ of $\mathcal{H}$ with $\xi(0,0)\in K$ satisfies
    \begin{align*}
        |\xi(t,j)|_{\mathcal{A}}\leq \beta(|\xi(0,0)|_{\mathcal{A}},t+j) + \kappa,
    \end{align*}
    for all $(t,j)\in\mathrm{dom}(\xi)$. \hfill $\QEDB$
\end{definition}

\section{Event-Triggered Hybrid Target Systems}
\begin{figure}
    \centering
    \includegraphics[width=0.95\linewidth]{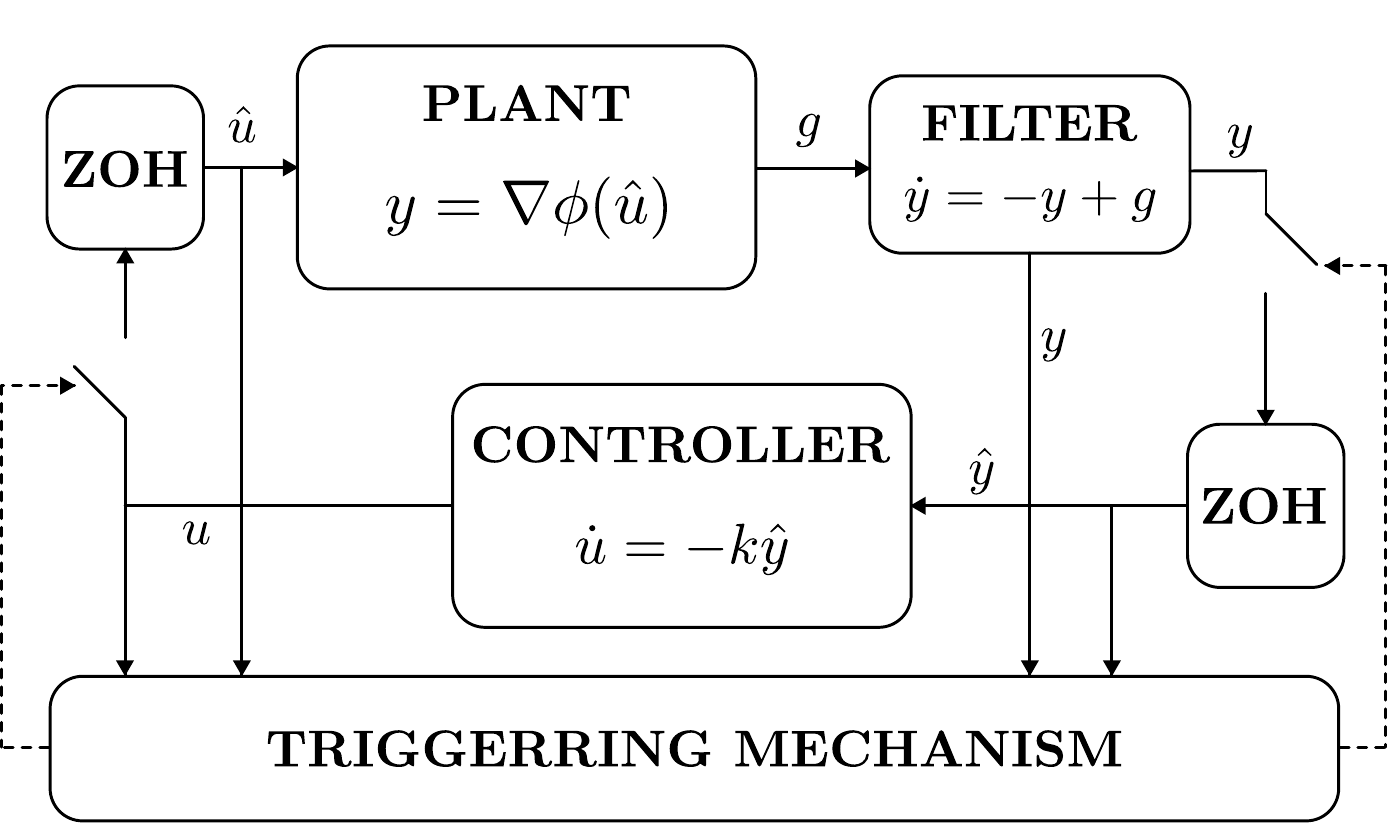}
    \caption{Block diagram of the event-triggered hybrid target system.}
    \label{fig:ET-GC-scheme}
\end{figure}
In this section, we introduce and study the behavior of the ``target system", i.e., the system that the event-triggered ES control scheme is designed to emulate. We call such a system an \emph{event-triggered hybrid target system}.
For simplicity, the plant output is taken to be the gradient of a (static) cost function $\phi$. In addition, throughout our manuscript, we assume that the cost function $\phi$ satisfies the following regularity assumption.
\begin{asmp}\label{asmp:cost-regularity}
    The function $\phi:\mathbb{R}^n\rightarrow\mathbb{R}$ is twice continuously differentiable, strongly convex, and there exists a constant $L_\phi>0$ such that
    \begin{align}
        \|\nabla\phi(u)-\nabla\phi(\hat{u})\|\leq L_\phi\|u-\hat{u}\|,
    \end{align}
    for all $u,\hat{u}\in\mathbb{R}^n$. \hfill $\QEDB$
\end{asmp}
The plant output $g=\nabla\phi(\hat{u})$ is passed to a low-pass filter with state $y$ that evolves according to 
\begin{align}
    \dot{y}=g-y.
\end{align}
The time-constant of the filter is taken to be equal to $1$ for simplicity but, through a suitable time-scale change, it can always be rescaled to $1$. Since the complexity of implementing a low-pass filter is not substantial, we treat it as part of the plant, although it can also be realized as a separate component attached to the physical plant.  
The output of the filter $y$ is passed to an event-triggering mechanism (defined below) which samples $y$ and forwards the sample to a zero-order hold (ZOH). The output $\hat{y}$ of the ZOH is taken to be the input to the controller which, for simplicity, is taken to be negative feedback with constant gain $k>0$, i.e., $u= - k \hat{y}$.
Then the output of the controller $u$ is also sampled by the event-triggering mechanism and forwarded to another ZOH whose output $\hat{u}$ serves as the plant input. In this way, the input to the plant $\hat{u}$ and the input to the controller $\hat{y}$ do not change, i.e., 
\begin{align}
    \dot{\hat{u}}&= 0, & \dot{\hat{y}}&= 0,
\end{align}
unless the event-triggering mechanism is activated, in which case they are updated according to the rules
\begin{align}
    \hat{u}^+&= u, & \hat{y}^+&= y.
\end{align}
A block diagram description of the closed-loop system described above is shown in Figure~\ref{fig:ET-GC-scheme}.
In order to study the closed-loop system, we formulate it as an HDS. We begin with the change of variables from $(u,y,\hat{u},\hat{y})$ to $(u,y,e_u,e_y)$ where 
\begin{align*}
    e_u&:=u-\hat{u}, & e_y&:=y-\hat{y}.
\end{align*}
We also define $x=(u,y)$ and $e=(e_u,e_y)$. We then define the HDS $\mathcal{H}=(C,F,D,G)$ with state $\xi=(x,e)$, and flow and jump maps
\begin{align*}
        F(\xi)&= \begin{pmatrix}f(x,e)\\ f(x,e)\end{pmatrix}, & G(\xi)&= \begin{pmatrix}x\\0\end{pmatrix},
    \end{align*}
    where the vector-valued functions $f$ is given by
    \begin{align}
        f(x,e):=\begin{pmatrix}- k(y-e_y)\\ \nabla\phi(u-e_u)-y\end{pmatrix}.
    \end{align}
The flow and jump maps of the HDS $\mathcal{H}$ capture the continuous evolution of $u$ and $y$, as well as the instantaneous updates of the samples $\hat{u}$ and $\hat{y}$. To complete the description of $\mathcal{H}$, we need to define the flow set $C$ and the jump set $D$, which will also define the event-triggering mechanism.

Throughout the remainder of the manuscript, the constant $\rho\in(0,\infty)$ serves as threshold that, if exceeded by the norm of the error $\|e\|$, activates the triggering mechanism. The corresponding flow and jump sets are defined by
\begin{subequations}\label{eq:first-mechanism}
    \begin{align}
        C&= \{\xi\in\mathbb{R}^{4n}~|~\|e\|^2\leq \rho\},\\
        D&= \{\xi\in\mathbb{R}^{4n}~|~\|e\|^2\geq \rho\}.
    \end{align}
\end{subequations}
Note that this triggering mechanism does not require time regularization (although time-triggered ES and periodic event-triggered ES can also be considered in our framework). We then have the following result, which is the cornerstone for the main results of this paper.

\vspace{0.1cm}
\begin{theorem}\label{thm:gc-is-ugas}
    There exists $k^*>0$ such that for any fixed $k\in(0,k^*)$, there exist constants $\beta_1,\beta_2>0$ such that, for all $\rho\in(0,\infty)$, and for all $\xi(0,0)\in C\cup D$, there exists a nontrivial solution, and every maximal solution is complete. Moreover, the compact set $\mathcal{A}(\rho)\subset C\cup D$ defined by
    \begin{align}\label{stableset}
        \mathcal{A}(\rho):=\{(x,e)~|~V(x)\leq 2\beta_1^{-1}\rho,\,\,\beta_2^2\|e\|^2\leq 4\beta_1 \rho\},
    \end{align}
    where $V$ is the function defined by
    \begin{align*}
        V(x):=(1-k)(\phi(u)-\phi(u^\star)) + \frac{1}{2}k \,\|y-\nabla\phi(u)\|^2,
    \end{align*}
    is UGAS for the HDS $\mathcal{H}$. \hfill $\QEDB$
\end{theorem}

\vspace{0.1cm}
\begin{rem}
    Our focus on the mechanism defined by \eqref{eq:first-mechanism} is only for illustration purposes and serves as a proof of concept. Indeed, more general event-triggering ES mechanisms, including dynamic ones, can be modeled in the framework of Hybrid Dynamical Systems adopted herein by leveraging existing results in the event-triggered control literature \cite{postoyan2014framework}. Future work will explore such generalizations. \hfill $\QEDB$
\end{rem}

%
Having established UGAS stability of the compact set $\mathcal{A}$ for the target HDS $\mathcal{H}$, we now move on to design an extremum-seeking algorithm that emulates the event-triggered hybrid target system $\mathcal{H}$. 
\section{Event-Triggered Extremum Seeking}
In the current section, we formulate and study the behavior of the proposed \emph{Event-Triggered ES} schemes. In particular, we consider two different schemes, shown in Figures 2 and 3, which differ in their tunable parameters and the analytical tools needed for their analysis.
\begin{figure}
    \centering
    \includegraphics[width=0.95\linewidth]{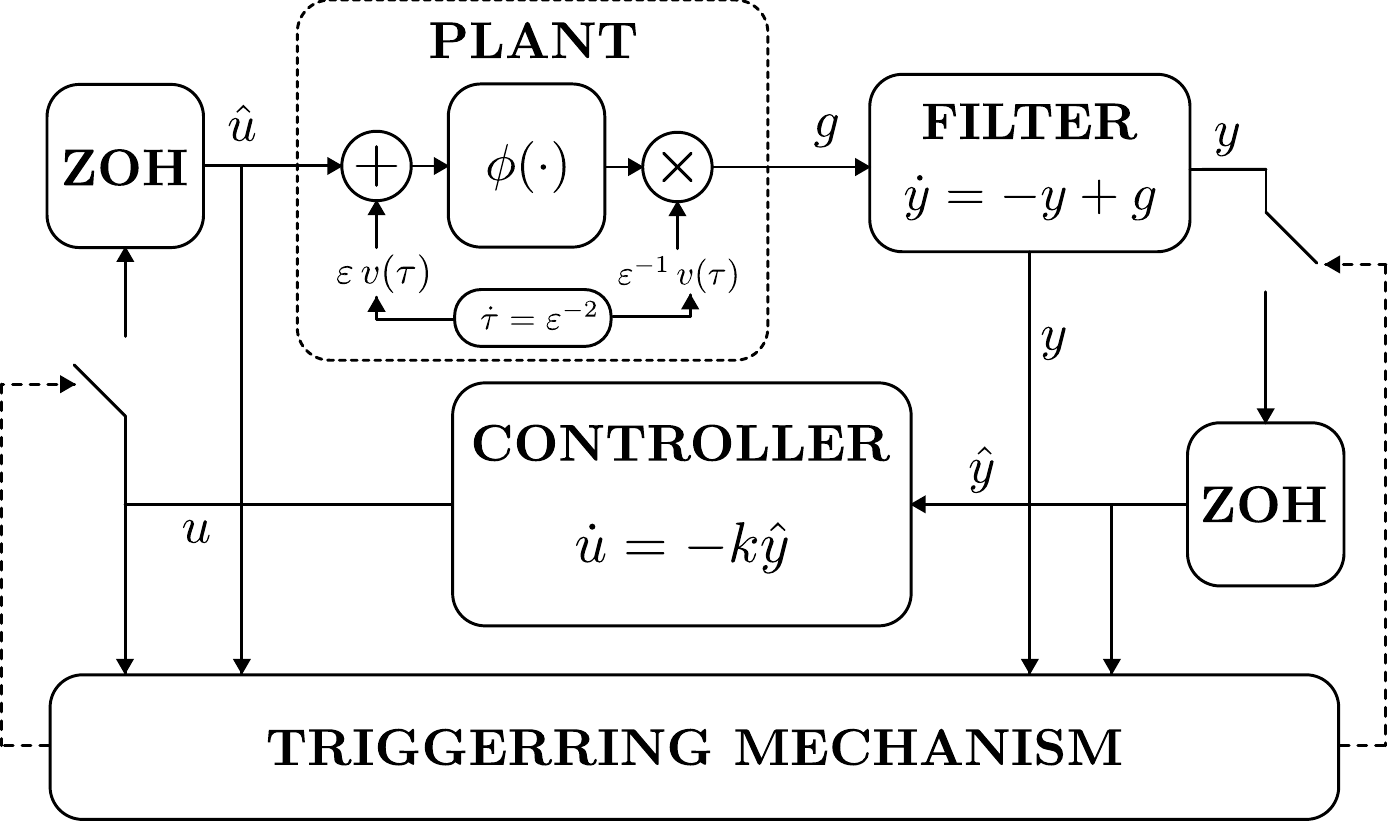}
    \caption{Block diagram of the event-triggered ES controller with tuning of the single parameter $\varepsilon$.}
    \label{fig:ET-ES-scheme}
\end{figure}
\subsection{Event-Triggered ES based on Lie-Bracket Averaging}
We begin with the block diagram description shown in Figure~\ref{fig:ET-ES-scheme}. In the diagram, the gradient of the static cost function $\nabla\phi$ is replaced with the dynamic term
\begin{align}
    g&= \frac{1}{\varepsilon}\phi(\hat{u}+\varepsilon v(\tau)) v(\tau), & \dot{\tau}&= \frac{1}{\varepsilon^2},
\end{align}
where $\varepsilon>0$ is a positive tuning parameter, and the function $v:\mathbb{R}\rightarrow\mathbb{R}^n$ is any continuous $T$-periodic vector-valued function satisfying
\begin{align}\label{periodicitydither}
    \int_0^T v(\tau)\,\mathrm{d}\tau&=0, & \frac{1}{T}\int_0^T v(\tau)v(\tau)^\top\,\mathrm{d}\tau&= I.
\end{align}
Note that $\varepsilon^{-2}$ describes the rate of change of $\tau$, which can be interpreted as the frequency in the classical ES schemes.

Using Hadamard's lemma \cite[Lemma 2.8]{nestruev2003smooth}, we obtain that
\begin{align*}
    g&= \frac{1}{\varepsilon} v(\tau) \phi(\hat{u})+ v(\tau)v(\tau)^\top\nabla \phi(\hat{u}) + R_{\varepsilon}(\hat{u},\tau),
\end{align*}
where the remainder term $R_{\varepsilon}$ is given explicitly by
\begin{align*}
    R_{\varepsilon}(\hat{u},\tau):=\int_0^1v(\tau)v(\tau)^\top (\nabla\phi(\hat{u}+\varepsilon \lambda  v(\tau))-\nabla\phi(\hat{u}))\,\mathrm{d}\lambda.
\end{align*}
The remaining components of the block diagram are the same as those in Figure~\ref{fig:ET-GC-scheme}. 
We now formulate the closed-loop system described by the block diagram in Figure~\ref{fig:ET-ES-scheme} as an HDS. Specifically, we define the HDS 
$$\mathcal{H}_{\varepsilon}=(\tilde{C},F_{\varepsilon},\tilde{D},\tilde{G}),$$
with the state $(\xi,\tau)$ where, as in the definition of the event-triggered hybrid target system HDS $\mathcal{H}$, the state $\xi$ is decomposed as $\xi=(x,e)$, with $x=(u,y)$, and $e=(e_u,e_y):=(u-\hat{u},y-\hat{y})$. The flow map of the HDS $\mathcal{H}_{\varepsilon}$ is defined by
\begin{align*}
    F_{\varepsilon}(\xi,\tau)&= \begin{pmatrix}\sum_{i=1}^2\varepsilon^{i-2} f_i(\xi,\tau) + f_\varepsilon(\xi,\tau),\varepsilon^{-2} \end{pmatrix},
\end{align*}
where the vector-valued functions $f_i$ and $f_{\varepsilon}$ are defined by
\begin{align*}
    f_1(\xi,\tau)&= \begin{pmatrix}
        0, v(\tau)\phi(u-e_u), 0, 0
    \end{pmatrix}, \\ 
    f_2(\xi,\tau)&= \begin{pmatrix}
        -k(y-e_y), v(\tau)v(\tau)^\top \nabla \phi(u-e_u)-y, 0,0
    \end{pmatrix},\\
    f_\varepsilon(\xi,\tau)&= \begin{pmatrix}
        0, R_\varepsilon(u-e_u,\tau), 0,0
    \end{pmatrix},
\end{align*}
and the jump map is defined by
\begin{align*}
    \tilde{G}(\xi,\tau)=(G(\xi),\tau),
\end{align*}
where $G$ is the jump map of the HDS $\mathcal{H}$, and the flow and jump sets 
\begin{align*}
    \tilde{C}&= C\times\mathbb{R}_{\geq 0}, & \tilde{D}&= D\times\mathbb{R}_{\geq 0},
\end{align*}
where $C$ and $D$ are the flow and jump sets of the HDS $\mathcal{H}$. 

The following theorem, which is the first main result of this paper, states the four key properties of the proposed hybrid ES controller: a) existence and completeness of solutions; b) absence of Zeno behavior; c) closeness of solutions with respect to the target system (on compact sets and compact time domains); and d) semi-global practical stability of the set of interest.

\vspace{0.1cm}
\begin{theorem}
    Let Assumption \ref{asmp:cost-regularity} be satisfied. For any $\rho>0$, define the set $\tilde{\mathcal{A}}(\rho)=\mathcal{A}(\rho)\times\mathbb{R}_{\geq 0}$, where $\mathcal{A}(\rho)$ is defined in \eqref{stableset}, and let $k^*>0$ be given by Theorem \ref{thm:gc-is-ugas}. Then, for any fixed $k\in(0,k^*)$, any $\rho>0$, any $T>0$, and any pair $\Delta>\nu>0$, there exists $\varepsilon^*>0$ such that for any $\varepsilon\in(0,\varepsilon^*)$ the following properties hold for every maximal solution $\zeta=(\xi,\tau)$ of $\mathcal{H}_{\varepsilon}$ with $|\zeta(0,0)|_{\tilde{\mathcal{A}}(\rho)}\leq\Delta$:
    \begin{enumerate}[(a)]
 \item \textsl{Existence and completeness of solutions:} $\zeta$ is well defined, and $\mathrm{dom}(\zeta)$ is unbounded, i.e.,  $\zeta$ is complete.
\item \textsl{Absence of Zeno behavior:} The bound \eqref{dwelltime} holds for $\zeta$, i.e., the solutions $\zeta$ to $\mathcal{H}_{\varepsilon}$ have a uniform semi-global dwell time.
\item \textsl{Closeness of solutions:} There exists a solution $\xi'$ to $\mathcal{H}$ such that the hybrid arcs $\xi'$ and $\xi$ are $(T,\epsilon)$-close.
\item \emph{Stability:} The solution $\zeta$ satisfies the following bound:
\begin{equation}
|\zeta(t,j)|_{\tilde{\mathcal{A}}(\rho)}\leq\beta(|\zeta(0,0)|_{\tilde{\mathcal{A}}(\rho)},t+j)+\nu,
\end{equation}
for all $(t,j)\in\mathrm{dom}(\zeta)$.\hfill \QEDB
\end{enumerate}
\end{theorem}

\vspace{0.1cm}

\begin{figure}
    \centering
    \includegraphics[width=0.95\linewidth]{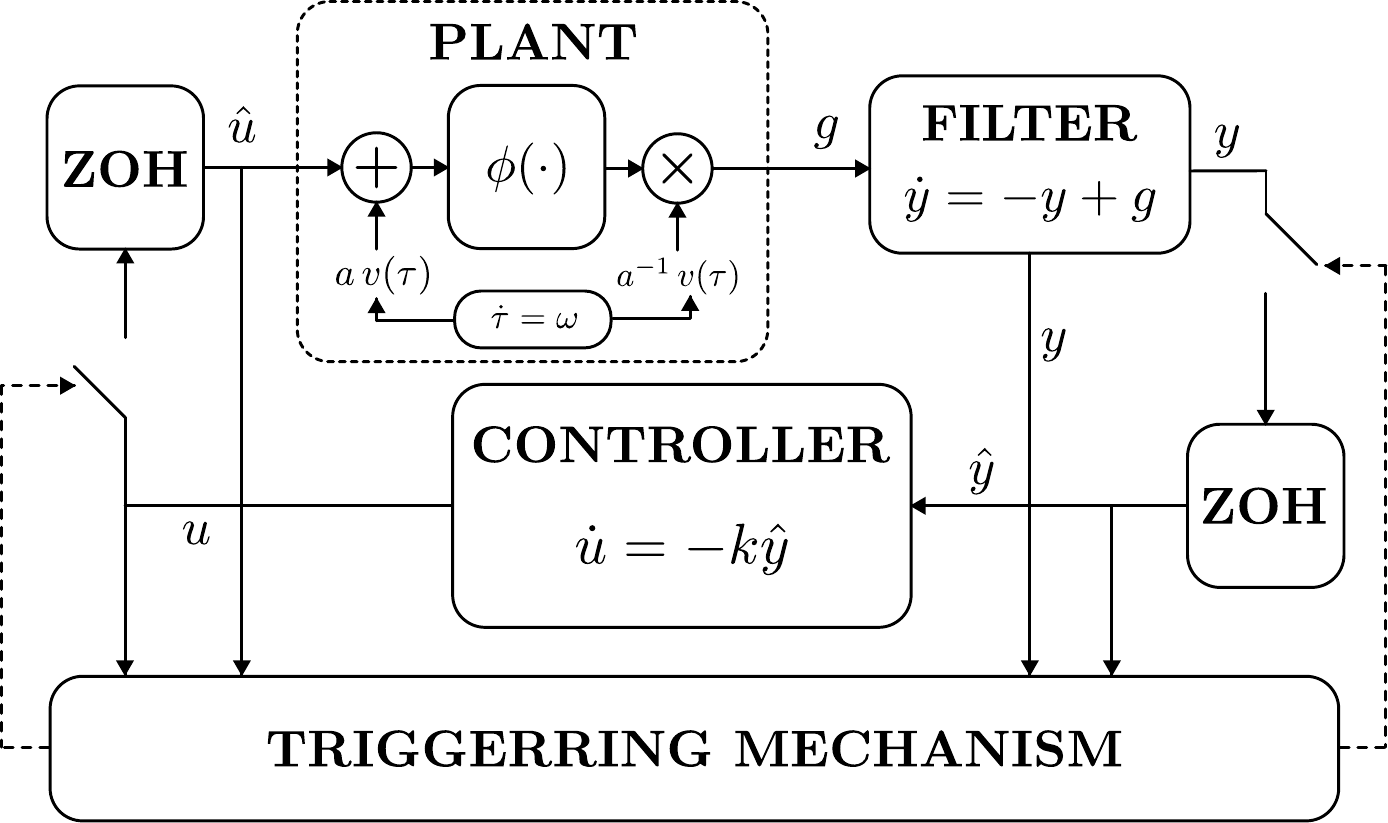}
    \caption{Block diagram of the event-triggered ES controller with sequential tuning of the parameters $a,\omega$.}
    \label{fig:ET-ES-scheme-classical}
\end{figure}

\subsection{Event-Triggered ES based on Standard Averaging}
It is also possible to study a variation of the proposed event-triggered ES using classical averaging by leveraging existing tools for hybrid ES \cite{poveda2017framework}. In particular, Figure~\ref{fig:ET-ES-scheme-classical} shows a block diagram that implements an ES where the gradient $\nabla\phi$ is obtained via the dynamics term:
\begin{align}
    g&= \frac{1}{a}\phi(\hat{u}+a v(\tau)) v(\tau), & \dot{\tau}&= \omega,
\end{align}
where $a>0$ and $\omega>0$ are now tunable parameters. Indeed, in a manner similar to the definition of the HDS $\mathcal{H}_{\varepsilon}$, the ES scheme in Figure~\ref{fig:ET-ES-scheme-classical} corresponds to the HDS $\mathcal{H}_{a,\omega}=(\tilde{C},F_{a,\omega},\tilde{D},\tilde{G}),$ 
where the flow and jump sets $C$ and $D$ and the jump map $G$ are the same as those of the HDS $\mathcal{H}_{\varepsilon}$ and the flow map $F_{a,\omega}$ is 
\begin{align*}
    F_{a,\omega}(\xi,\tau)=\begin{pmatrix} f_{a,\omega}(\xi,\tau),\omega\end{pmatrix},
\end{align*}
and the vector-valued map $f_{a,\omega}$ is given explicitly by
\begin{align*}
    f_{a,\omega}(\xi,\tau)&=\begin{pmatrix}
        0\\ a^{-1}v(\tau)\phi(u-e_u)\\ 0\\0
    \end{pmatrix} \\
    &+\begin{pmatrix}
        -k(y-e_y)\\ v(\tau) v(\tau)^\top \nabla \phi(u-e_u)-y\\ 0\\0
    \end{pmatrix} + \mathcal{O}(a).
\end{align*}
In this case, the proposed scheme can be studied using the framework of \cite{poveda2017framework} by noticing that the dither $v$ satisfying \eqref{periodicitydither} can be equivalently generated via a time-invariant dynamic oscillator with gain $\omega$, as in \cite[Eq. (7)]{poveda2017framework}. Using the results from \cite[Sec.7.2]{poveda2017framework} the average flow and jump maps of the event-triggered ES system are now given by:
\begin{align*}
    F_a(\xi)&= \begin{pmatrix}-k (y-e_y),\nabla \phi(u-e_u)-y+\mathcal{O}(a),0, 0\end{pmatrix}, \\
    G_a(\xi)&= \begin{pmatrix}u,y,0,0\end{pmatrix},
\end{align*}
\begin{figure*}[!htb]
\centering
\begin{subfigure}{0.325\textwidth}
\centering
\includegraphics[width=1\linewidth]{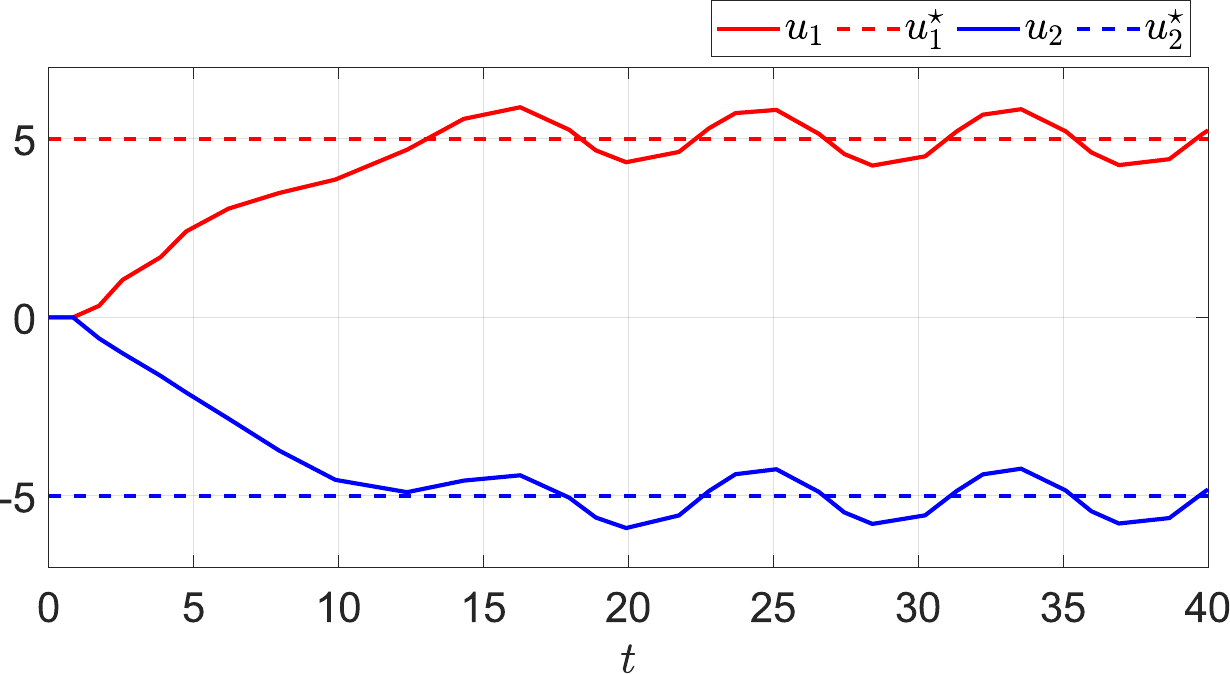} 
\caption{}
\label{fig:subim1}
\end{subfigure}
\begin{subfigure}{0.325\textwidth}
\centering
\includegraphics[width=1\linewidth]{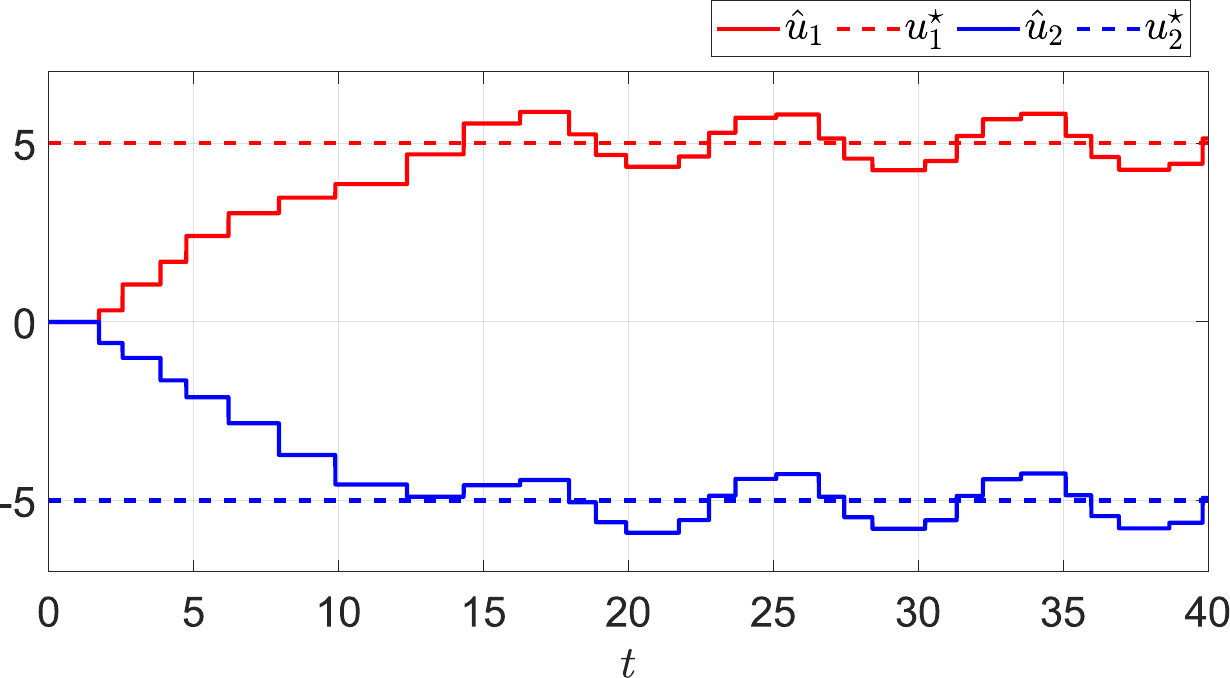}
\caption{}
\label{fig:subim2}
\end{subfigure}
\begin{subfigure}{0.325\textwidth}
\centering
\includegraphics[width=1\linewidth]{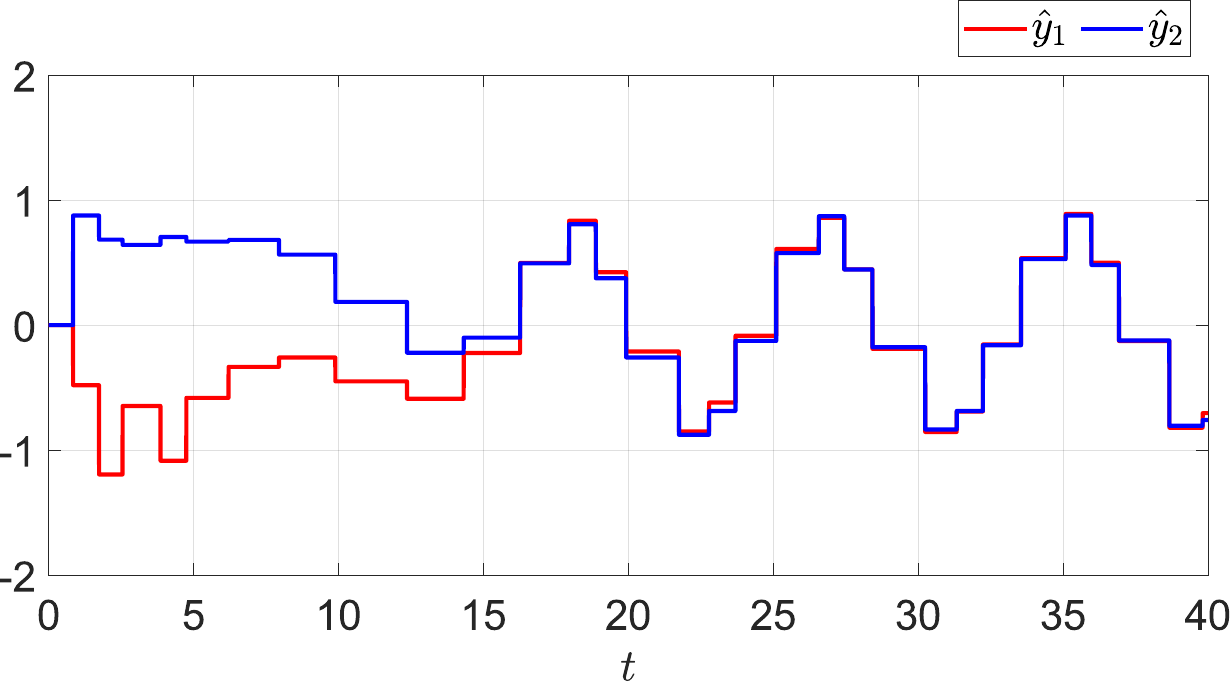}
\caption{}
\label{fig:subim4}
\end{subfigure}
\caption{Numerical simulation results for the proposed event-triggered ES.}
\label{fig:example-simulations}
\end{figure*}
and the flow and jump sets remain the same as in \eqref{eq:first-mechanism}. The following theorem, which is our second main result, formalizes the properties of the HDS $\mathcal{H}_{a,\omega}$. 

\vspace{0.1cm}
\begin{theorem}
 Let Assumption \ref{asmp:cost-regularity} be satisfied. For any $\rho>0$, define the set $\tilde{\mathcal{A}}(\rho)=\mathcal{A}(\rho)\times\mathbb{R}_{\geq 0}$, where $\mathcal{A}(\rho)$ is defined in \eqref{stableset}, and let $k^*>0$ be given by Theorem \ref{thm:gc-is-ugas}. Then, for any fixed $k\in(0,k^*)$, any $\rho>0$, $T>0$, and each pair $\Delta>\nu>0$ there exists $a^*>0$ such that for all $a\in(0,a^*)$ there exists $\omega^*>0$ such that for all $\omega>\omega^*$ and all solutions $\zeta=(\xi,\tau)$ to $\mathcal{H}_{a,\omega}$ with $|\zeta(0,0)|_{\tilde{\mathcal{A}}(\rho)}\leq\Delta$, the following four properties hold:
 \begin{enumerate}[(a)]
\item \textsl{Existence of solutions:} $\zeta$ is well defined, and $\mathrm{dom}(\zeta)$ is unbounded, i.e.,  $\zeta$ is complete.
\item \textsl{Absence of Zeno behavior:} The bound \eqref{dwelltime} holds for $\zeta$, i.e., the solutions $\zeta$ to $\mathcal{H}_{a,\omega}$ have a uniform semi-global dwell time.
\item \textsl{Closeness of solutions:} There exists a solution $\xi'$ to $\mathcal{H}$ such that $\xi'$ and $\xi$ are $(T,\epsilon)$-close.
\item \emph{Stability:} The solution $\zeta$ satisfies the following bound:
\begin{equation}
|\zeta(t,j)|_{\tilde{\mathcal{A}}(\rho)}\leq\beta(|\zeta(0,0)|_{\tilde{\mathcal{A}}(\rho)},t+j)+\nu,
\end{equation}
for all $(t,j)\in\mathrm{dom}(\zeta)$.\hfill \QEDB
\end{enumerate}
\end{theorem}

\vspace{0.1cm}
In the next section, we illustrate via simulations the properties of the proposed event-triggered ES controllers.

\section{Numerical Simulations}
Having proven the stability of the proposed event-triggered ES scheme under the natural conditions of Assumption \ref{asmp:cost-regularity}, we now provide numerical simulations that validate our results. 

Throughout this section, we take $n=2$ and consider the cost function defined by
\begin{subequations}\label{eq:cost_function_example_1}
\begin{align}
    \phi(u)&=\frac{1}{2} (u-u^\star)^\top Q (u-u^\star),
\end{align}
where the vector $u^\star$ and the symmetric matrix $Q$ are 
\begin{align}
    u^\star&= \begin{pmatrix} 5 \\ -5 \end{pmatrix}, &  Q&=\frac{1}{2}\begin{pmatrix} \sqrt{2} & 1\\ 1 & \sqrt{2}\end{pmatrix}.
\end{align}
\end{subequations}
Clearly, the function $\phi$ satisfies Assumption \ref{asmp:cost-regularity}. For the parameters defining the HDS $\mathcal{H}$, we take 
\begin{align}
    \rho&=1, & k&= 0.75.
\end{align}
We take the vector-valued function $v$ to be
\begin{align*}
    v(\tau)=\sqrt{2}(\sin(2\pi \tau),\sin(4\pi\tau)),
\end{align*}
and use the parameter $\varepsilon=0.04$ in the HDS $\mathcal{H}_{\varepsilon}$. Numerical simulation results are shown in Figure~\ref{fig:example-simulations}. As demonstrated by the results figures, the samples $\hat{u}$ and $\hat{y}$ are updated only occasionally as intended. As expected, the solutions of the system do not converge to the equilibrium point $(u^\star,0,u^\star,0)$ precisely due to the error induced by the sampling. Nevertheless, it is clear that the parameters converge to a neighborhood of said point. Moreover, the size of this neighborhood can be shrank by reducing the positive constant $\rho$ at the cost of higher computational cost (i.e., higher sampling rate). 

\section{Conclusion}
We introduced a class of event-triggered hybrid extremum-seeking controllers. The algorithms are naturally modeled as hybrid dynamical systems, for which suitable regularity, stability, and convergence properties can be derived using well-established tools in the literature. The main challenge in the design of the controller is to obtain suitable triggering rules that do not depend on the high frequency of the controller, while still leading to a well-posed hybrid system. The results were presented for static mappings but can be extended in the future to incorporate dynamic plants, as well as other types of triggering rules and/or optimization dynamics modeled as well-posed hybrid dynamical systems with suitable regularity and stability properties, as discussed in \cite{poveda2017framework} and \cite{abdelgalil2025lie}.

\balance 
\bibliographystyle{IEEEtran}
\bibliography{reference.bib}

\end{document}